\documentclass[letterpaper,11pt]{amsart}

\usepackage[margin=1.2in]{geometry}
\usepackage{amsmath,amsthm,amssymb, mathtools}
\usepackage{xspace,xcolor}
\usepackage[breaklinks,colorlinks,citecolor=teal,linkcolor=teal,urlcolor=teal,pagebackref,hyperindex]{hyperref}
\usepackage[alphabetic]{amsrefs}
\usepackage{comment}
\usepackage[all]{xy}
\usepackage{tikz-cd}
\usepackage{color}
\usepackage{comment}

\theoremstyle{plain}
\newtheorem{theo}{Theorem}[section]

\newtheorem{lemm}[theo]{Lemma}
\newtheorem{prop}[theo]{Proposition}

\newtheorem{conj}{Conjecture}

\theoremstyle{definition}
\newtheorem{defi}[theo]{Definition}

\theoremstyle{remark}
\newtheorem{rema}[theo]{Remark}

\newcommand{\bfA}{\mathbf{A}}

\newcommand{\bb}[1]{\mathbb{#1}}

\newcommand{\PP}{\mathbb{P}}
\newcommand{\QQ}{\mathbb{Q}}

\newcommand{\RR}{\mathbb{R}}
\newcommand{\CC}{\mathbb{C}}

\newcommand{\VV}{\mathbb{V}}
\newcommand{\ZZ}{\mathbb{Z}}

\newcommand{\cD}{\mathcal{D}}

\newcommand{\cH}{\mathcal{H}}

\newcommand{\cL}{\mathcal{L}}

\newcommand{\cO}{\mathcal{O}}

\newcommand{\uind}[1]{^{(#1)}}

\newcommand{\ubind}[1]{^{[#1]}}

\newcommand{\lsta}{_{*}}

\newcommand{\sta}{^{*}}

\newcommand{\norm}[1]{\|#1\|}

\newcommand{\uround}[1]{\left\lceil #1 \right\rceil }
\newcommand{\lround}[1]{\left\lfloor #1 \right\rfloor}

\DeclareMathOperator{\rank}{rank}

\DeclareMathOperator{\codim}{codim}

\DeclareMathOperator{\Sym}{Sym}

\DeclareMathOperator{\im}{im}

\DeclareMathOperator{\reg}{reg}

\DeclareMathOperator{\supp}{supp}

\DeclareMathOperator{\Baily}{BB}

\usepackage{soul}

\begin{document}
	
	\subjclass[2020]{14C30, 14D07, 14E30, 14J27}
	
	\author{Hyunsuk Kim}
	
	\address{Department of Mathematics, University of Michigan, 530 Church Street, Ann Arbor, MI 48109, USA}
	
	\email{kimhysuk@umich.edu}
	
	\begin{abstract} 
	     Fujino gave a proof in \cite{Fujino-CBF} for the semi-ampleness of the moduli part in the canonical bundle formula in the case when the general fibers are K3 surfaces or Abelian varieties. We show a similar statement when the general fibers are primitive symplectic varieties with mild singularities. This answers a question of Fujino raised in the same article. Moreover, using the structure theory of varieties with trivial first Chern class, we reduce the question of semi-ampleness in the case of families of K-trivial varieties to a question when the general fibers satisfy a slightly weaker Calabi-Yau condition.
	\end{abstract}
	
	\title[A Remark on Fujino's work on the canonical bundle formula]{A Remark on Fujino's work on the canonical bundle formula via period maps}
	
	\maketitle
	
 

        \section{Introduction}
        Starting from Kodaira's canonical bundle formula for minimal elliptic surfaces in \cite{Kodaira}, the canonical bundle formula has been extensively studied and widely generalized (for example, by \cite{Kawamata-codim2}, \cite{Kawamata-subadjunction-general}, \cite{Ambro04-Boundary}, \cite{FM-canonical-bundle-formula}, \cite{Fujino-CBF}, \cite{Filipazzi-gen-adjunction}, to list a few) and crucially used in higher dimensional adjunction. Roughly speaking, the set up is the following. Let $f \colon X \to Y$ be a projective morphism with connected fibers between normal projective varieties. Let $\Delta$ be a (non necessarily effective) $\QQ$-divisor on $X$ such that $(X, \Delta)$ is klt (or lc) and $K_{X} + \Delta \sim_{\QQ, f} 0$. Then the canonical bundle formula suggests to write
        $$ K_{X} + \Delta \sim_{\QQ} f\sta (K_{Y} + B_{Y} + M_{Y})$$
        in an insightful way so that, roughly speaking, $B_{Y}$ comes from the contribution of the singularity of the fibers, and $M_{Y}$ comes from the (Hodge theoretic) variation of the general fiber. Moreover, once we have $f \colon X\to Y$, it makes sense to have the same formula for other birational models $f' \colon X' \to Y'$ over $f$ so that $\mathbf{B}_{Y}$ and $\mathbf{M}_{Y}$ make sense as $b$-divisors. We call $\mathbf{B}_{Y}$ the divisorial part and $\mathbf{M}_{Y}$ the moduli part. After choosing a high enough model satisfying certain simple normal crossing assumptions, the moduli part commutes with generically finite base change and it is known to be nef (\cite{Kawamata-subadjunction-general}). In the following cases for the general fibers, this nef divisor is known to be semi-ample:
        \begin{enumerate}
            \item $\PP^{1}$: \cite{Kawamata-codim2}, \cite{Pro-Sho:semiampleness}*{Theorem 8.1}
            \item Elliptic curves: \cite{Kodaira}, \cite{Fujita-elliptic}
            \item Surfaces with $\kappa(X_{\eta}) = 0$ or Abelian varieties: \cite{Fujino-CBF}
            \item Surfaces not isomorphic to $\PP^{2}$: \cite{Filipazzi-gen-adjunction}
            \item $\PP^{2}$, and hence in relative dimension 2: \cite{ascher2023moduli}*{Theorem 1.4}.
        \end{enumerate}
        A conjecture of Prokhorov and Shokurov \cite{Pro-Sho:semiampleness}*{Conjecture 7.13} predicts that this would always be the case. They suggest to compactify the space parametrizing the general fibers to a projective variety and show that the ample line bundle on the compactified moduli space pulls back to the moduli part.
        
        Meanwhile, \cite{Filipazzi-gen-adjunction} gives an inductive approach to this conjecture by developing the canonical bundle formula and adjunction for generalized pairs. Using the techniques in the minimal model program, he shows that in order to verify the conjecture in relative dimension $n$, it is enough to verify the conjecture in relative dimension $< n$, and the following two extremal cases in relative dimension $n$:
        \begin{enumerate}
            \item $f\colon X \to Y$ is a Mori fiber space, or
            \item $K_{X} \sim_{\QQ, f} 0$, and $\Delta$ has no horizontal divisors.
        \end{enumerate}

        We concentrate on special situations in the second extremal case, when $\Delta$ has no horizontal divisors. Indeed, the cases when the general fibers are K3 or Abelian varieties in \cite{Fujino-CBF} fit into this framework. Fujino observes that the semi-ampleness in this situation reduces to a purely Hodge theoretic result, namely, by looking at the Baily-Borel compactification of the period domain, we get semi-ampleness. He raises the question whether one can use the same strategy when the fiber is a holomorphic symplectic variety. We provide an affirmative answer.

        \begin{theo} \label{theo:hyperkahler-semiample}
            Let $f \colon X \to Y$ be an algebraic fiber space, with $X$ and $Y$ projective normal varieties. Suppose that the general fiber is a primitive symplectic variety of dimension $2m$ (see Definition \ref{defi:prim-Symplectic}). Then there exists a vertical divisor $\Delta$ on $X$ such that $f : (X, \Delta) \to Y$ is an lc-trivial fibration. The moduli part $\mathbf{M}_{Y}$ does not depend on $\Delta$ and it is $b$-semi-ample.
        \end{theo}

        We use the same idea as in \cite{Fujino-CBF}, namely, we consider the Baily-Borel compactification of the parametrizing space of weight 2 Hodge structures. The other input is to relate the second cohomology and the middle cohomology of the general fibers. While it is hard to determine the entire middle cohomology of the fiber, the existence of a generically non-degenerate two form is strong enough for our purpose since we are only interested in the variation of the lowest piece of the Hodge filtration.
        
        We can combine this result with the structure theory of klt varieties with trivial first Chern class, and reduce the semi-ampleness question to the case when the general fibers satisfy a weaker Calabi-Yau condition. First, we recall the statement of the conjecture of Prokhorov and Shokurov.

        \begin{conj}[\cite{Pro-Sho:semiampleness}*{Conjecture 7.13}]
            Let $f: (X, \Delta) \to Y$ be an lc-trivial fibration. Then
            \begin{enumerate}
                \item $\mathbf{M}_{Y}$ is $b$-semi-ample.
                \item Let $X_{\eta}$ be the generic fiber of $f$. Then there exists an $I_{0}$ depending only on the dimension of the generic fiber and the coefficients of the horizontal part $\Delta^{h}$ such that $I_{0} (K_{X_{\eta}} + \Delta_{\eta}) \sim 0$.
                \item $\mathbf{M}_{Y}$ is effectively $b$-semi-ample, that is, there exists an $I_{1}$ depending only on the dimension of the generic fiber and the coefficients of $\Delta^{h}$ such that $I_{1}\mathbf{M}_{Y}$ is $b$-free.
            \end{enumerate}
        \end{conj}
        
        We present two subconjectures of the first part of \cite{Pro-Sho:semiampleness}*{Conjecture 7.13} in our setting.

        \begin{conj} \label{conj:semi-ample-CY}
            Let $f \colon X \to Y$ be an algebraic fiber space, with $X$ and $Y$ normal projective varieties. Consider a pair $(X, \Delta)$ such that $K_{X} + \Delta \sim_{\QQ, f} 0$ and $\Delta$ has no horizontal divisors. Suppose that the general fiber of $f$
            \begin{enumerate}
                \item is a Calabi-Yau manifold, or
                \item is a pre-CY variety (see Definition \ref{defi:singular-HK-CY}).
            \end{enumerate}
            Then $ \mathbf{M}_{Y}$ is $b$-semi-ample.
        \end{conj}

        \begin{conj} \label{conj:semi-ample-allK-triv}
            Let $f \colon X \to Y$ be an algebraic fiber space as before. Consider the pair $(X, \Delta)$ such that $K_{X} + \Delta \sim_{\QQ, f} 0$ and $\Delta$ has no horizontal divisors.
            Suppose that the general fiber
            \begin{enumerate}
                \item is smooth, or
                \item has klt singularities.
            \end{enumerate}
            Then $\mathbf{M}_{Y}$ is $b$-semi-ample.
        \end{conj}

        We show the following result:
        \begin{theo} \label{theo:reduction-to-CY}
            Conjecture \ref{conj:semi-ample-CY} (1) (resp. (2)) in reletive dimension $\leq n$ implies Conjecture \ref{conj:semi-ample-allK-triv} (1) (resp. (2)) in relative dimension $n$.
        \end{theo}

        Therefore, we reduce the $b$-semi-ampleness question for $K$-trivial fibrations to the (pre)-CY case.

        \section{Preliminaries}
        
        \subsection{The Canonical bundle formula} We discuss the canonical bundle formula and the behavior of the moduli part under various operations, mainly following \cite{Kollar-KodCBF}. Before that, we collect some notation and terminology for algebraic fiber spaces.

        {\bf Notation and terminology.}
        \begin{enumerate}
            \item An {\it algebraic fiber space} is a projective morphism $f\colon X\to Y$ between normal projective varieties with connected fibers. We put $\dim f = \dim X - \dim Y$ for the {\it relative dimension} of the algebraic fiber space.
            
            \item For a $\QQ$-divisor $B$ on a normal algebraic variety $X$, we write $B = B_{+}-B_{-}$ where $B_{+}$ and $B_{-}$ share no irreducible components. We call $B_{+}$ (resp. $B_{-}$) the {\it positive} (resp. {\it negative}) part of $B$.
            
            \item Let $f \colon X\to Y$ be an algebraic fiber space and $B$ be a $\QQ$-divisor on $X$. We write $B = B^{h} + B^{v}$ where the irreducible components of $B^{h}$ are exactly the irreducible components of $B$ that map onto $Y$. We call $B^{h}$ (resp. $B^{v}$) the {\it horizontal} (resp. {\it vertical}) part of $B$. We say $B$ is {\it horizontal} (resp. {\it vertical}) if $B = B^{h}$ (resp. $B = B^{v}$).

            \item \cite{Ambro04-Boundary} We say $f \colon (X, \Delta) \to Y$ is an {\it lc-trivial fibration} if $f \colon X \to Y$ is an algebraic fiber space such that
            \begin{enumerate}
                \item $K_{X} + \Delta \sim_{\QQ, f} 0$,
                \item $(X, \Delta)$ is klt over the generic point of $Y$, and
                \item If $\pi \colon (X', \Delta') \to (X, \Delta)$ is a log resolution, then $h^{0}(F', \uround{\Delta_{-}'|_{F}}) = 1$, where $F'$ is a general fiber of $f \circ \pi \colon X' \to Y$. \footnote{This is equivalent to the condition $\rank f\lsta \cO_{X}(\uround{\bfA (X, \Delta)}) = 1$ in \cite{Ambro04-Boundary}*{Definition 2.1}.}
            \end{enumerate}
             
            \item Let $f \colon X \to Y$ be an algebraic fiber space and let $\mu \colon Y' \to Y$ be a projective morphism from a normal projective variety $Y'$. Let $X'$ be any normal projective variety mapping birationally onto the main component $X' \to (X \times_{Y} Y')_{\mathrm{main}}$. Then we have the corresponding commutative diagram
            $$ \begin{tikzcd}
                X' \ar[d, "f'"] \ar[r, "\mu'"] & X \ar[d, "f"] \\ Y' \ar[r , "\mu"] & Y.
            \end{tikzcd}$$
            We call $f' \colon X' \to Y'$ an {\it algebraic fiber space induced by} $\mu$.

            \item Let $f\colon X\to Y$ be an algebraic fiber space and $\Delta$ be a divisor on $X$ such that $K_{X} + \Delta \sim_{\QQ, f} 0$. We say that $f$  {\it satisfies the standard normal crossing assumptions} if the following conditions hold:
            \begin{enumerate}
                \item $X$ and $Y$ are smooth,
                \item There exists an SNC divisor $\Sigma$ on $Y$ such that $f$ is smooth over $Y \setminus \Sigma$,
                \item $\supp(\Delta) + f\sta \Sigma$ has SNC support,
                \item $\Delta$ is relatively SNC over $Y \setminus \Sigma$.
            \end{enumerate}
            Even though the divisors $\Delta$ and $\Sigma$ are oppressed in the terminology, we remark that they are the part of the data of the standard normal crossing assumptions. For clarity, we will sometimes say that $X, Y, \Delta, \Sigma$ satisfy the standard normal crossing assumptions.

            \item Let $f \colon X \to Y$ be an algebraic fiber space satisfying the standard normal crossing assumptions. We say $f$ is {\it semi-stable in codimension 1} if there exists a codimension $\geq 2$ closed subset $Z$ of $Y$ such that $f\sta (\Sigma \setminus Z)$ is a reduced SNC divisor. In this case, the local system $R^{i}f\lsta \CC_{X}|_{Y \setminus \Sigma}$ has unipotent monodromy for every $i$ by \cite{Katz-regularity}.
        \end{enumerate}

        \begin{rema}
            For any algebraic fiber space $f\colon X \to Y$, there exists a birational morphism $\mu \colon Y' \to Y$ and an algebraic fiber space $f' \colon X' \to Y'$ induced by $\mu$ such that $f'$ is satisfies the standard normal crossing assumptions. Moreover, we can take a generically finite morphism $\mu\colon Y' \to Y$ and an algebraic fiber space $f' \colon X' \to Y'$ induced by $\mu$ such that $f'$ is semi-stable in codimension 1.
        \end{rema}

        \begin{rema}
            Even though there is a more general set up such as slc-trivial fibrations due to \cite{Fujino-slc-trivial} which is important in many applications, we are only interested here in the case of lc-trivial fibrations.
        \end{rema}

        We describe the canonical bundle formula for algebraic fiber spaces. Let $f \colon (X, \Delta) \to Y$ be an lc-trivial fibration and fix $L$ a $\QQ$-divisor on $Y$ such that $K_{X} +\Delta \sim_{\QQ} f\sta L$. For each prime divisor $P$ on $Y$, we consider
        $$ t_{P} = \sup \{ t \in \QQ : (X, \Delta+ t f\sta P) \text{ is log-canonical over } \eta_{P} \}.$$
        We define $B_{Y} = \sum_{P} (1 - t_{P}) P$. Note that this is a finite sum since $(X, \Delta)$ is klt over the generic point of $Y$. We define $M_{Y} = L - K_{Y} - B_{Y}$, so that we have
        $$ K_{X} + \Delta \sim_{\QQ} f\sta (K_{Y} + B_{Y} + M_{Y}).$$
        
        \begin{rema} \label{rema:moduli-upto-res-of-X}
            Note that by definition, $B_{Y}$ and $M_{Y}$ do not depend on the birational model of $X$. In other words, if we consider a birational morphism $\pi \colon \widetilde{X} \to X$ such that $K_{\widetilde{X}} + \widetilde{\Delta} = \pi\sta (K_{X} + \Delta)$, then $f \colon (\widetilde{X} , \widetilde{\Delta}) \to Y$ is an lc-trivial fibration, and the $B_{Y}$ and $M_{Y}$ computed in terms of $f\circ \pi \colon (\widetilde{X}, \widetilde{\Delta}) \to Y$ agree with the ones computed in terms of $f \colon (X, \Delta) \to Y$. 
        \end{rema}
        
        For every birational morphism $\mu \colon Y' \to Y$, we have a commutative diagram 
        $$ \begin{tikzcd}
            X' \ar[r, "\mu'"] \ar[d, "f'"] & X \ar[d, "f"] \\ Y' \ar[r,"\mu"] & Y
        \end{tikzcd}$$
        where $\mu' \colon X' \to X$ is birational. If we write $K_{X'} + \Delta' = {\mu'}\sta (K_{X} + \Delta)$, then we use the formula for $f' \colon (X', \Delta') \to Y'$ and write
        $$ K_{X'} + \Delta' \sim_{\QQ} {f'}\sta(\mu\sta L) = {f'}\sta (K_{Y'} + B_{Y'} + M_{Y'}).$$
        We have $\mu\lsta M_{Y'} = M_{Y}$ and $\mu\lsta B_{Y'} = B_{Y}$ and therefore, we may and will consider the divisorial part and the moduli part as $b$-divisors, and denote them by $\mathbf{B}_{Y}$ and $\mathbf{M}_{Y}$. We point out that each $M_{Y'}$ may be well-defined only up to $\QQ$-linear equivalence, but once we fix a representative $L$ such that $K_{X} + \Delta \sim_{\QQ} f\sta L$, then $\mathbf{M}_{Y}$ is well-defined as a $b$-divisor.

        \begin{rema} \label{rema:moduli-std-behavior}
            We collect some standard facts about the divisorial part and the moduli part from \cite{Kollar-KodCBF}*{\S 8.4.}. We say that a $b$-divisor $\mathbf{D}$ on $X$ {\it descends to} $X'$ if $\mathbf{D} = \overline{\mathbf{D}_{X'}}$. We remark that a $b$-divisor $\mathbf{D}$ on $X$ is {\it $b$-nef} (resp. {\it $b$-semi-ample, $b$-free}) if there exists a birational model $X' \to X$ such that $\mathbf{D}$ descends to $X'$ and $\mathbf{D}_{X'}$ is nef (resp. semi-ample, free). For basic notions for $b$-divisors, we refer to \cite{flips-3fold4fold}*{Chapter 1}.
        \begin{enumerate}
            \item For an lc trivial fibration $f:(X, \Delta) \to Y$, we can take a resolution of singularities $\mu : Y' \to Y$ and an algebraic fiber space $f' :  (X', \Delta') \to Y'$ induced by $\mu$ such that $X'$, $Y'$, $\Delta'$, $B_{Y'}$ satisfy the standard normal crossing assumptions. In this case, $\mathbf{M}_{Y}$ and $\mathbf{K}+ \mathbf{B}_{Y}$ descend to $Y'$ in the sense that for any birational morphism $\pi \colon Y'' \to Y'$, we have
            $$ \mathbf{M}_{Y,Y''} = \pi\sta M_{Y'}, \quad \text{and} \quad K_{Y''} + \mathbf{B}_{Y, Y''} = \pi\sta (K_{Y'} + B_{Y'}).$$
            \item $M_{Y}$ only depends on the general fiber $(F, \Delta|_{F})$ and $Y$.
            \item  Let $f \colon (X, \Delta) \to Y$ be an lc-trivial fibration satisfying the standard simple normal crossing assumptions. Then $\mathbf{M}_{Y}$ is $b$-nef.
            \item Let $f \colon (X, \Delta) \to Y$ be an lc-trivial fibration satisfying the standard normal crossing assumptions. Let $\mu \colon Y' \to Y$ be a generically finite morphism from a smooth variety $Y'$. Let
            $$ \begin{tikzcd}
                X' \ar[r , "\mu'"] \ar[d, "f'"] & X \ar[d, "f"] \\ Y' \ar[r, "\mu"] & Y
            \end{tikzcd}$$
            be an algebraic fiber space induced by $\mu$ and write $K_{X'} + \Delta' = {\mu'}\sta (K_{X} + \Delta)$. Then, $f' \colon(X', \Delta') \to Y'$ is an lc-trivial fibration and
            $$ M_{Y'} = \mu\sta M_{Y}. $$
            \item Suppose moreover that there exists an SNC divisor $\Sigma$ on $ Y$ such that $f \colon (X, \Delta) \to Y$ satisfies the standard normal crossing assumptions and $\Delta^{h}$ is an integral divisor (hence $- \Delta^{h}$ is effective). Furthermore, assume that $p_{g}(X_{\eta}) = 1$ where $\eta$ is the generic point of $Y$, and the variation of Hodge structures $R^{\dim f}f\lsta \CC_{X}|_{Y \setminus \Sigma}$ has unipotent local monodromies. Then $M_{Y}$ is the divisor class corresponding to the canonical extension of the lowest piece of the Hodge filtration of $R^{\dim f}f\lsta \CC_{X}|_{Y \setminus \Sigma}$.
        \end{enumerate} 
        \end{rema}

        \begin{rema} \label{rema:check-semi-ample-gen-fini}
            We remark that it is enough to take a resolution of singularities and a generically finite cover of the base in order to check the $b$-semi-ampleness of the moduli part.
        \end{rema}

        We compare our setup with the formulation in \cite{FM-canonical-bundle-formula} and \cite{Fujino-CBF}. Consider an algebraic fiber space $f \colon X\to Y$ between smooth projective varieties. Suppose that the Kodaira dimension of the generic fiber of $f$ is zero, that is, $\kappa(X_{\eta}) =0$, where $\eta$ is the generic point of $Y$. Fix $b \in \ZZ_{>0}$ such that the $b$-th plurigenus of the general fiber $P_{b}(X_{\eta})$ is non-zero. Then we have the following formula for the canonical bundle $K_{X}$.

         \begin{prop}[\cite{FM-canonical-bundle-formula}*{Proposition 2.2.}] \label{prop:CBF-FM}
            In the above situation, there exists a unique $\QQ$-divisor $D$ on $Y$, modulo linear equivalence, with an isomorphism of graded $\cO_{Y}$-algebras:
            $$ \bigoplus_{i \geq 0} \cO_{Y} (\lround{iD}) \simeq \bigoplus_{i \geq 0} (f\lsta \cO_{Y}(ibK_{X/ Y}))^{\ast\ast}.$$
            Furthermore, the isomorphism induces a $\QQ$-linear equivalence
            $$ bK_{X} \sim_{\QQ} f\sta (bK_{Y} + D) + B,$$
            where $B$ is a $\QQ$-divisor on $X$ satisfying
            \begin{enumerate}
                \item $f\lsta \cO_{X}(\lround{iB_{+}}) = \cO_{Y}$ for $i > 0$, and
                \item $\codim_{Y} f(\supp B^{-}) \geq 2$.
            \end{enumerate}
        \end{prop}

        Note that we recover the canonical bundle formula in this situation, since $f\colon(X, -b^{-1} B) \to Y$ is an lc-trivial fibration, and we can write
        $$ K_{X} - b^{-1} B \sim_{\QQ} f\sta (K_{Y} + B_{Y} + M_{Y}).$$
           
        \begin{rema}
            We point out that in this situation, there is a natural choice of $\Delta^{h}$ since on the general fiber $F$, the divisor $-b\Delta|_{F}$ is the zero locus of the unique section (up to scalar) in $H^{0}(F, \omega_{F}^{\otimes b})$. Moreover, any two such $\Delta$ that make $(X,\Delta) \to Y$ as an lc-trivial fibration differ by the pull-back of a divisor on $Y$. Note that the moduli part $M_{Y}$ only depends on the general fiber $(F, \Delta|_{F})$. Therefore, given a projective morphism $f \colon X \to Y$ between smooth projective varieties whose general fiber has Kodaira dimension zero, it makes sense to talk about the moduli part $M_{Y}$ {\bf without} picking a divisor $\Delta$ on $X$ such that $f\colon(X, \Delta) \to Y$ is an lc-trivial fibration.
        \end{rema}

        We end this section by describing the behavior of the moduli part after taking a generically finite cover of the source.
        \begin{prop}[\cite{Fujino-CBF}*{Lemma 4.1}] \label{prop:Fujino-cover-formula}
            Let $f\colon X \to Y$ and $h \colon W\to Y$ be algebraic fiber spaces between smooth projective varieties such that
            \begin{enumerate}
                \item $\kappa(X_{\eta}) = 0$, where $\eta$ is the generic point of $Y$,
                \item there is a generically finite morphism $g \colon W \to X$ such that $h = f\circ g$,
                \item there is an SNC divisor $\Sigma$ on $Y$ such that $f$ and $h$ are smooth over $Y^{\circ} := Y \setminus \Sigma$, and
                \item $\kappa(W_{\eta}) =0$ and $p_{g}(W_{\eta}) = 1$.
            \end{enumerate}
            Let $M_{X/Y}$ and $M_{W/Y}$ be the moduli part of the canonical bundle formula coming from $f\colon X \to  Y$ and $h \colon W \to Y$ respectively. Then $M_{X/ Y} = M_{W/ Y}$.
        \end{prop}

        \begin{rema}
        We point out that there is a difference of a multiple of $b$ in the formula from \cite{Fujino-CBF}, where $b$ is the smallest number such that the plurigenus $P_{b}(X_{\eta})$ is non-zero. This is because the semi-stable part denoted by $L_{X/Y}^{ss}$ in \cite{Fujino-CBF} actually equals $b M_{X/Y}$ in our situation.
        \end{rema}
        
    \subsection{The Structure theorem for $K$-trivial klt varieties}
    We recall the structure theorem for $K$-trivial varieties, starting from the decomposition theorem of Beauville-Bogomolov and its singular generalization. Roughly speaking, this theorem suggests to study $K$-trivial varieties by studying three different special types of varieties.
    
    \begin{theo} \cite{Beauville-decomposition}
        Let $X$ be a Kähler manifold with $c_{1}(K_{X}) = 0 \in H^{2}(X, \RR)$. Then $X$ admits a finite étale cover $\gamma \colon \widetilde{X} \to X$ such that $\widetilde{X}$ decomposes as
        $$ \widetilde{X} \simeq A \times \prod_{j \in J} Y_{j} \times \prod_{k \in K} Z_{k} $$
        such that
        \begin{enumerate}
            \item $A$ is an Abelian variety;
            \item $Y_{j}$ are irreducible hyperkähler manifolds;
            \item $Z_{k}$ are Calabi-Yau manifolds.
        \end{enumerate}
    \end{theo}

    We introduce a remarkable generalization of this result to singular varieties, due to \cite{HP:klt-beauville-bogo-decomposition} and a series of works including \cite{Greb-Guenancia-Kebekus:kltBeauville-Bog}, \cite{Druel-dectheorem}, \cite{Greb-Kebekus-Peternell:singspaceswtrivcan}. We first define the singular analogues of irreducible hyperkähler manifolds and Calabi-Yau manifolds, in the sense of \cite{Greb-Guenancia-Kebekus:kltBeauville-Bog}. We also introduce weaker versions of these notions.

    \begin{defi} \label{defi:singular-HK-CY}
        Let $X$ be a normal projective variety of dimension $n \geq 2$. We say
        \begin{enumerate}
            \item $X$ is {\it CY (Calabi-Yau)} if $X$ has Gorenstein canonical singularities with $\omega_{X} \simeq \cO_{X}$, and if $H^{0}(Y, \Omega_{Y}\ubind{p}) = 0$ for all covers $\gamma \colon Y \to X$ which are étale in codimension 1 and for all $1 \leq p \leq n -1$.
            \item[(1)'] $X$ is {\it pre-CY (pre-Calabi-Yau)} if $X$ has Gorenstein canonical singularities with $\omega_{X} \simeq \cO_{X}$, and if $H^{0}(X, \Omega_{X}\ubind{p}) = 0$ for all $1 \leq p \leq n-1$.
            \item $X$ is {\it IHS (irreducible holomorphic symplectic)} if $X$ has Gorenstein canonical singularities with $\omega_{X} \simeq \cO_{X}$, and if there exists a holomorphic 2-form $\sigma \in H^{0}(X, \Omega_{X}\ubind{2})$ such that for all covers $\gamma \colon Y \to X$ étale in codimension 1, the exterior algebra $H^{0}(Y, \Omega_{Y}\ubind{\bullet})$ is generated by the reflexive pull-back of $\sigma$.
        \end{enumerate}
    \end{defi}

     Following \cite{Beau-Symsing}, we also give a slightly general class of symplectic varieties than those appearing in the decomposition theorem.

    \begin{defi}\label{defi:prim-Symplectic}
        A normal projective variety $X$ is {\it primitive symplectic} if 
        \begin{enumerate}
            \item $H^{1}(X, \cO_{X}) =0$ and $H^{0}(X, \Omega_{X}\ubind{2}) = \CC \sigma$, where $\sigma$ is non-degenerate on the smooth locus $X_{\reg}$, and
            \item there exists a resolution of singularities $\pi \colon Y \to X$ such that the the pull back of $\sigma|_{X_{\reg}}$ extends to a holomorphic 2-form on $Y$.
        \end{enumerate}
    \end{defi}

    \begin{rema} \label{rema:relate-H2-and-middle}
        Let $X$ be a primitive symplectic variety of dimension $2m$ and $\sigma$ be the symplectic form. Then $\omega_{X} \simeq \cO_{X}$ since it is trivialized by $\sigma^{m}$ and therefore $X$ is Gorenstein. We point out that the second condition on Definition \ref{defi:prim-Symplectic} is equivalent to $X$ having canonical singularities by \cite{GKKP:extension-of-holomorphic-forms}*{Theorem 1.4}. Pick any resolution of singularities $\pi \colon \widetilde{X} \to X$ and let $\widetilde{\sigma}$ be the holomorphic $2$-form on $\widetilde{X}$ extending $\pi\sta\sigma|_{X_{\reg}}$. Note that in this case, we have
        $$ h^{2,0} (\widetilde{X}) = h^{0,2}(\widetilde{X}) = 1.$$
        Note that there is a natural morphism of Hodge structures $\mu \colon \Sym^{m}H^{2}(\widetilde{X}, \QQ) \to H^{2m}(\widetilde{X}, \QQ)$. We observe that $\ker \mu$ is a Hodge structure of weight $2m$ that does not have a $(2m, 0)$ part since $\widetilde{\sigma}$ is generically non-degenerate. We also point out that if $X$ itself is a smooth hyperkähler manifold, then the natural morphism $\mu$ for $X$ is injective by \cite{Verb96}.
    \end{rema}

    \begin{rema}
        We will use the following convention for Calabi-Yau manifolds. A smooth projective variety $X$ of dimension $n$ is a Calabi-Yau manifold if $\omega_{X} \simeq \cO_{X}$ and
        $$ \dim_{\CC} H^{k,0}(X) = \begin{cases}
            1 & \text{if } k = 0, n \\ 0 & \text{otherwise.}
        \end{cases}$$
        This condition is weaker than the condition appearing in Beauville's decomposition theorem (see \cite{Beauville-decomposition}*{Proposition 2}). However, we point out that if $X$ is a holomorphic symplectic manifold of dimension $2m$ with
        $$ \dim_{\CC} H^{k,0}(X) = \begin{cases}
            1 & k \text{ is even} \\ 0 & k \text{ is odd,}
        \end{cases}$$
        then $X$ is automatically simply connected, and hence irreducible holomorphic symplectic (see \cite{Hnw}*{Proposition A.1}).
    \end{rema}

    We finally introduce the generalization of the Beauville-Bogomolov decomposition to the singular case.

    \begin{theo}[\cite{HP:klt-beauville-bogo-decomposition}*{Theorem 1.5}] \label{HP:singular-decomposition-theorem}
        Let $X$ be a normal projective variety with at worst klt singularities such that $c_{1}(K_{X}) = 0$. Then there exists a projective variety $X'$ with at worst canonical singularities, with a quasi-étale (which means quasi-finite and étale in codimension one) map $\gamma\colon X' \to X$ and a decomposition
        $$ X' \simeq A \times \prod_{j \in J} Y_{j} \times \prod_{k\in K} Z_{k},$$
        into normal varieties with trivial canonical bundles, such that
        \begin{enumerate}
            \item $A$ is an Abelian variety;
            \item $Y_{j}$ are irreducible holomorphic symplectic varieties;
            \item $Z_{k}$ are Calabi-Yau varieties.
        \end{enumerate}
    \end{theo}

    \subsection{Summary of Fujino's result for K3 surfaces}
        We briefly summarize the Hodge theoretic results used in \cite{Fujino-CBF}*{\S 2}. Let $B$ be a smooth projective variety and $\Sigma$ be an SNC divisor on $B$. Let $B^{\circ} = B \setminus \Sigma$ and consider a polarized $\ZZ$-variation of Hodge structures of weight 2 on $B^{\circ}$ with the following numerical conditions:
        $$ h^{2,0} = h^{0,2} = 1, \qquad h^{1,1} =g \geq 3, \qquad \text{and } h^{p,q} = 0 \text{ otherwise}.$$
        We furthermore assume that the local monodromies around $\Sigma$ are unipotent and there exists a neat arithmetic group $\Gamma$ containing the local monodromy operators around $\Sigma$. Then we have the period map $\wp^{\circ}\colon B^{\circ} \to \cD/ \Gamma$ and in this case, $\cD$ is a bounded hermitian symmetric domain. By Borel's extension theorem \cite{Borel-extension}, the holomorphic map $\wp^{\circ}$ extends to $B$ as
        $$ \wp \colon B \to (\cD / \Gamma)^{\Baily}$$
        where $(\cD/ \Gamma)^{\Baily}$ is the Baily-Borel compactification of $\cD/ \Gamma$ (\cite{Baily-Borel-compactification}) which is a normal analytic space. The tautological sub-bundle on $\cD$ descends to a line bundle $\cL$ on $\cD/ \Gamma$. The sections of $\cL^{\otimes gk}$ can be identified with automorphic forms of weight $k$ (which are $\Gamma$-equivariant $k$-pluricanonical forms on $\cD$). For some $k > 0$, the automorphic forms give an embedding of $\cD/ \Gamma$ in a projective space, and the automorphic forms can be continuously (hence analytically) extended to $(\cD/ \Gamma)^{\Baily}$. Moreover, these extended automorphic forms define a projective embedding of $(\cD/ \Gamma)^{\Baily}$. Hence, $\cL^{\otimes gk}$ extends to an ample line bundle $\cO_{(\cD/ \Gamma)^{\Baily}}(1)$ on $(\cD/\Gamma)^{\Baily}$. On the other hand, from the variation of Hodge structures, we have the associated vector bundle $\cH$ on $B^{\circ}$ with a filtration $F^{\bullet}$. By the nilpotent orbit theorem, the canonical extension $\overline{\cH}$ of $\cH$ is a vector bundle on $B$ which carries a filtration $F^{\bullet}$ by vector bundles extending the filtration on $\cH$. By definition, we have the natural identification $(\wp^{\circ})\sta \cL \simeq F^{2}\cH$. \cite{Fujino-CBF}*{Theorem 2.10} says that the lowest piece of the filtration on the canonical extension $\overline{\cH}$ and the ample line bundle on $(\cD/ \Gamma)^{\Baily}$ are compatible. In other words, we have
        $$ \wp\sta \cO_{(\cD/\Gamma)^{\Baily}}(1) \simeq (F^{2}\overline{\cH})^{\otimes gk}.$$
        In particular, this shows that $F^{2} \overline{\cH}$ is a semi-ample line bundle on $B$.

        \begin{rema}
            \cite{Fujino-CBF}*{\S 2} has a parallel statement dealing with variations of Hodge structures of weight 1 which covers the case when the general fiber is an Abelian variety.
        \end{rema}

        \section{Proof of Theorem \ref{theo:hyperkahler-semiample}}
        We give a proof of Theorem \ref{theo:hyperkahler-semiample}.

        \begin{proof}[Proof of Theorem \ref{theo:hyperkahler-semiample}]
            All other assertions except $b$-semi-ampleness are clear. By Remark \ref{rema:check-semi-ample-gen-fini}, we can take a generically finite base change of $Y$ and resolve singularities. Hence, we can assume that we have a morphism $f'\colon X' \to Y'$ between smooth projective varieties and we have a divisor $\Delta'$ on $X'$ such that $f' \colon (X', \Delta') \to Y'$ is an lc-trivial fibration. We can furthermore assume that the followings hold:
            \begin{enumerate}
                \item there exists an SNC divisor $\Sigma'$ on $Y'$ such that $f'$ is smooth over $Y' \setminus \Sigma'$, and
                \item $X', Y', \Delta', \Sigma'$ satisfy the standard normal crossing assumptions, and $f'$ is semi-simple in codimension 1, and
                \item every fiber $X_{y}'$ of $f'$ for $y \in Y' \setminus \Sigma$ admits a unique holomorphic 2-form $\sigma_{y} \in H^{0}(X_{y}', \Omega_{X_{y}'}^{2})$ which is generically non-degenerate, and
                \item $h^{1,1}(X_{y}') \geq 3$ (by blowing-up $X'$ further).\footnote{We point out that it is important to allow negative coefficients for ${\Delta'}^{h}$ in the canonical bundle formula since we are blowing-up $X'$ further.}
            \end{enumerate}
            Note that in this case, the local monodromies of $R^{2}f'\lsta \CC_{X'}|_{Y'\setminus \Sigma'}$ and $R^{2m}f'\lsta \CC_{X'}|_{Y' \setminus \Sigma'}$ around $\Sigma'$ are unipotent. It is enough to show that $M_{Y'}$ is semi-ample. Since the general fibers of $f \colon X \to Y$ have Gorenstein canonical singularities, $-{\Delta'}^{h}$ is effective and integral. Hence, we are exactly in the situation in Remark \ref{rema:moduli-std-behavior} (5). Hence, $M_{Y'}$ is the divisor class of the canonical extension of the lowest piece of the Hodge filtration of the variation of Hodge structures $R^{2m}f'\lsta \CC_{X'}|_{Y' \setminus \Sigma'}$. We denote by $\VV\uind{2}$ and $\VV\uind{2m}$ the polarizable variations of Hodge structures $R^{2}f'\lsta \CC_{X'}|_{Y' \setminus \Sigma'}$ and $R^{2m}f'\lsta \CC_{X'}|_{Y' \setminus \Sigma'}$.\footnote{We point out that we should first fix a relatively ample class for $f' :X' \to Y'$ to give a polarization on the variations of Hodge structures $\VV\uind{2}$ and $\VV\uind{2m}$ due to a subtlety of signs in the polarization for primitive and non-primitive parts, but this does not cause a problem throughout the argument.} Then we have a natural morphism of variations of Hodge structures
            $$ \mu: \Sym^{m} \VV\uind{2} \to \VV\uind{2m}.$$
            Then we have non-canonical splittings
            \begin{align*}
                \Sym^{m} \VV\uind{2} &\simeq \ker \mu \oplus \im \mu \\
                \VV\uind{2m} & \simeq \im \mu \oplus \bb{B}
            \end{align*}
            where $\ker \mu, \im \mu,$ and $\bb{B}$ are polarizable variations of Hodge structures. By Remark \ref{rema:relate-H2-and-middle}, $\ker \mu$ and $\bb{B}$ do not have $(2m, 0)$-part. We denote by $(\cH\uind{2}, F^{\bullet})$, $(\cH, F^{\bullet})$, and $(\cH\uind{2m}, F^{\bullet})$ the filtered vector bundles associated to the variations of Hodge structures $\VV\uind{2}, \im\mu,$ and $\VV\uind{2m}$. Denote by $\overline{\cH\uind{2}}, \overline{\cH},$ and $ \overline{\cH\uind{2m}}$ the canonical extensions of $\cH\uind{2}, \cH,$ and $ \cH\uind{2m}$, respectively. Since $\bb{B}$ does not have a $(2m, 0)$-part, we get
            $$ F^{2m} \overline{\cH} \simeq F^{2m} \overline{\cH\uind{2m}}.$$
            We also have $F^{2m}\overline{\cH} \simeq (F^{2}\overline{\cH\uind{2}})^{\otimes m}$ by Lemma \ref{lemm:comparison-H2-Sym} below and since $\ker \mu$ does not have a $(2m, 0)$-part. Then \cite{Fujino-CBF}*{Theorem 2.10} immediately implies that $F^{2}\overline{\cH\uind{2}}$ is semi-ample, and therefore, $F^{2m}\overline{\cH\uind{2m}}$ is semi-ample as well. Since $M_{Y'}$ is a divisor class of this line bundle, we are done.
        \end{proof}

        \begin{lemm} \label{lemm:comparison-H2-Sym}
            Let $Y$ be a smooth complex manifold and $\Sigma$ be an SNC divisor on $Y$. Let $\VV$ be a polarizable variation of Hodge structures of weight 2 on $Y \setminus \Sigma$, with Hodge numbers
            $$ h^{2,0} = h^{0,2} =1, \quad \text{and} \quad h^{p,q} = 0 \text{ if $p < 0$ or $p > 2$.} $$
            Suppose that $\VV$ has unipotent local monodromies along $\Sigma$. Let $\cH$ and $\cH'$ be the filtered vector bundles associated to the variations of Hodge structures $\VV$ and $\Sym^{m}\VV$, respectively. Denote by $\overline{\cH}$ and $\overline{\cH'}$ be the canonical extensions of $\cH$ and $\cH'$, respectively. Then we have
            $$ (F^{2}\overline{\cH})^{\otimes m} \simeq F^{2m} \overline{\cH'}.$$
        \end{lemm}

        \begin{proof}
            We fix a polarization on $\VV$ and denote the hermitian metric on $\cH$ by $h$. We have the induced metric $h'$ on $\cH'$. Note that there is a canonical isomorphism $(F^{2} \cH)^{\otimes m} \simeq F^{2m} \cH'$ between line bundles on $Y \setminus \Sigma$. Fix local coordinates $z_{1},\ldots, z_{n}$ on an open subset $\Omega \subset Y$ such that $\Omega \simeq \{ (z_{1},\ldots, z_{n}) \in \CC^{n} : |z_{i}| < 1\}$ and $\Sigma \cap \Omega$ is given by the equation $z_{1}\cdots z_{l} = 0$. Since $\VV$ and $\Sym^{m}\VV$ both have unipotent local monodromies, we have the following description of the local sections of $F^{2} \overline{\cH}$ and $F^{2m} \overline{\cH'}$ using the estimation of Hodge norms by \cite{Schmid-VHS}*{\S 6}:
            \begin{align*}
                \Gamma(\Omega, F^{2} \overline{\cH}) & = \left\{ s\in \Gamma(\Omega \setminus \Sigma, F^{2}\cH) : \norm{s}_{h}^{2} = O\left(\prod_{i=1}^{l}|\log z_{i}|^{N_{i}} \right) \text{for some $N_{i} > 0$} \right\} \\
                \Gamma(\Omega, F^{2m} \overline{\cH'}) & = \left\{ s' \in \Gamma(\Omega \setminus \Sigma, F^{2m}\cH') : \norm{s'}_{h'}^{2} = O\left(\prod_{i=1}^{l}|\log z_{i}|^{N_{i}'} \right) \text{for some $N_{i}' > 0$} \right\}.
            \end{align*}
            Since $\norm{s^{\otimes m}}_{h'}^{2} = \norm{s}_{h}^{2m}$, the isomorphism $(F^{2} \cH)^{\otimes m} \simeq F^{2m} \cH'$ extends to $(F^{2} \overline{\cH})^{\otimes m} \simeq F^{2m} \overline{\cH'}$.
        \end{proof}

        \section{Proof of Theorem \ref{theo:reduction-to-CY}}
        We give a proof of Theorem \ref{theo:reduction-to-CY}.

        \begin{proof}[Proof of Theorem \ref{theo:reduction-to-CY}]
            We first consider the case when the general fiber is smooth. By resolving singularities, we can assume that $X$ and $Y$ are both smooth, and the geometric general fiber has trivial first Chern class. By \cite{Beauville-decomposition}, after taking a generically finite base change of $Y$, we can assume that there are algebraic fiber spaces $h \colon W \to Y$, $p \colon A \to Y, \{q_{j} \colon Y_{j} \to Y\}_{j =1}^{a}, \{r_{k}\colon Z_{k} \to Y\}_{k=1}^{b}$ satisfying the following conditions:
            \begin{enumerate}
                \item $h$ factors as $ W \xrightarrow{g} X \xrightarrow{f} Y$, where $g$ is a generically finite morphism.
                \item There exists an SNC divisor $\Sigma$ on $Y$ such that $f, h, p, q_{j}, r_{k}$ are smooth over $Y^{\circ} = Y \setminus \Sigma$.
                \item Over $Y^{\circ}$, the fiber spaces $p, q_{j}, r_{k}$'s are families of Abelian varieties, irreducibe hyperkähler manifolds, and Calabi-Yau manifolds, respectively.
                \item Denote the inverse images of $Y^{\circ}$ by $W^{\circ}$, $A^{\circ}$, $Y_{j}^{\circ}$, and $Z_{k}^{\circ}$. Then there is an isomorphism
                $$ W^{\circ} \simeq A^{\circ} \times_{Y^{\circ}} Y_{1}^{\circ} \times_{Y^{\circ}} \ldots \times_{Y^{\circ}} Y_{a}^{\circ} \times_{Y^{\circ}} Z_{1}^{\circ} \times_{Y^{\circ}} \ldots \times_{Y^{\circ}} Z_{b}^{\circ} $$
                over $Y^{\circ}$.
                \item The variations of Hodge structures $R^{\dim h}h\lsta \CC_{W^{\circ}}, R^{\dim p}p\lsta \CC_{A^{\circ}}, R^{\dim q_{j}}q_{j\ast} \CC_{Y_{j}^{\circ}}$, and $ R^{\dim r_{k}} r_{k \ast} \CC_{Z_{k}^{\circ}}$ have unipotent local monodromies around $\Sigma$.
            \end{enumerate}
            By Proposition \ref{prop:Fujino-cover-formula}, it is enough to show that the moduli part $M_{W/ Y}$ associated to $h \colon W \to Y$ is semi-ample. Applying the Künneth formula and the realization of the moduli part as a canonical extension of the lowest piece of the Hodge filtration, we get
            $$ M_{W/Y}= M_{A/ Y} + \sum_{j=1}^{a} M_{Y_{j}/Y} + \sum_{k=1}^{b} M_{Z_{k}/ Y},$$
            where $M_{A/ Y}, M_{Y_{j}/ Y},$ and $M_{Z_{k}}$ are the moduli parts associated to the fiber spaces $p, q_{j}$, and $r_{k}$, respectively. By \cite{Fujino-CBF}*{Theorem 5.1} and Theorem \ref{theo:hyperkahler-semiample}, we know that $M_{A/ Y}$ and $M_{Y_{j}/ Y}$ are semi-ample. Hence, Conjecture \ref{conj:semi-ample-CY} (1) in relative dimension $ \leq n$ implies Conjecture \ref{conj:semi-ample-allK-triv} (1) in relative dimension $n$.

            We deal with the singular case similarly. By applying Theorem \ref{HP:singular-decomposition-theorem} to the geometric generic fiber $X_{\overline{\eta}}$, we get the following diagram:
            $$ \begin{tikzcd}
                & & & \widetilde{X}  \ar[d, "\pi"]  \\
                X'''  \ar[r, "\mu_{1}"] \ar[rrru, bend left = 15, "h"] & X''  \ar[r, "\mu_{2}"] & X'  \ar[r, "\mu_{3}"] & X_{\overline{\eta}} 
            \end{tikzcd} $$
            such that
            \begin{enumerate}
                \item $\pi \colon \widetilde{X} \to X_{\overline{\eta}} $ is a resolution of singularities,
                \item $\mu_{3}$ is a generically finite morphism such that
                $$ X' \simeq A \times \prod_{j=1}^{a} Y_{j} \times \prod_{k=1}^{b} Z_{k},  $$
                where $A$ is an Abelian variety, $Y_{j}$ are IHS varieties, and $Z_{k}$ are CY varieties,
                \item $\phi_{j} : \widetilde{Y}_{j} \to Y_{j}$ and $\varphi_{k}: \widetilde{Z}_{k} \to Z_{k}$ are resolutions of singularities and
                $$ X'' \simeq A \times \prod_{j=1}^{a} \widetilde{Y}_{j} \times \prod_{k=1}^{b} \widetilde{Z}_{k}$$
                such that $X'' \to X'$ is induced by the maps $\phi_{j}$ and $\varphi_{k}$'s.
                \item $X'''$ is smooth and $\mu_{1} : X''' \to X''$ is a birational morphism resolving the indeterminacy of the rational map from $X''$ to $\widetilde{X}$.
            \end{enumerate}
            
            After replacing $Y$ with a generically finite cover of $Y$ and resolving singularities, we can assume that $Y$ is smooth, and there exist an SNC divisor $\Sigma$ and algebraic fiber spaces $p \colon A \to Y$, $\{ q_{j} \colon Y_{j} \to Y\}_{j=1}^{a}$, $\{ \widetilde{q}_{j} 
            \colon Y_{j} \to Y\}_{j=1}^{a}$, $\{ r_{k} \colon Z_{k}\to Y\}_{k=1}^{b}$, and $\{ \widetilde{r}_{k} \colon \widetilde{Z}_{k}\to Y\}_{k=1}^{b}$ over $Y$, and a commutative diagram
            $$ \begin{tikzcd}
                & & & \widetilde{X} \ar[rd, "\widetilde{f}"] \ar[d, "\pi"] & \\
                X''' \ar[r, "\mu_{1}"] \ar[rrru, bend left = 15, "h"] & X'' \ar[r, "\mu_{2}"] & X' \ar[r, "\mu_{3}"] & X \ar[r, "f"] & Y
            \end{tikzcd} $$
            satisfying the following conditions:
            \begin{enumerate}
                \item $\pi \colon \widetilde{X} \to X$ is a resolution of singularities, and $\widetilde{f}$ satisfies the standard normal crossing assumptions.
                \item For every point $y \in Y \setminus \Sigma$, the fibers of $p$, $q_{j}$, and $r_{k}$ over $y$ are Abelian varieties, primitive symplectic varieties, and pre-CY varieties, respectively.
                \item $\mu_{3} \colon X' \to X$ is a generically finite morphism and if we denote by ${X'}^{\circ}$, $A^{\circ}, Y_{j}^{\circ}, Z_{k}^{\circ}$ the corresponding inverse images of $Y^{\circ}$, then we have an isomorphism
                $$ {X'}^{\circ} \simeq A^{\circ} \times_{Y^{\circ}} Y_{1}^{\circ} \times_{Y^{\circ}} \ldots \times_{Y^{\circ}} Y_{a}^{\circ} \times_{Y^{\circ}} Z_{1}^{\circ} \times_{Y^{\circ}} \ldots \times_{Y^{\circ}} Z_{b}^{\circ} $$
                over $Y^{\circ}$.
                \item $\mu_{2}$ is a birational morphism from a smooth projective variety $X''$, and if we denote the inverse images of $Y^{\circ}$ similarly, we have an isomorphism and a commutative diagram
                $$ \begin{tikzcd}
                    {X''}^{\circ} \ar[r ,"\simeq"] \ar[d,"\mu_{2}"] & A^{\circ} \times_{Y^{\circ}} \widetilde{Y}_{1}^{\circ} \times_{Y^{\circ}} \ldots \times_{Y^{\circ}} \widetilde{Y}_{a}^{\circ} \times_{Y^{\circ}} \widetilde{Z}_{1}^{\circ} \times_{Y^{\circ}} \ldots \times_{Y^{\circ}} \widetilde{Z}_{b}^{\circ} \ar[d] \\
                    {X'}^{\circ} \ar[r, "\simeq"] & A^{\circ} \times_{Y^{\circ}} Y_{1}^{\circ} \times_{Y^{\circ}} \ldots \times_{Y^{\circ}} Y_{a}^{\circ} \times_{Y^{\circ}} Z_{1}^{\circ} \times_{Y^{\circ}} \ldots \times_{Y^{\circ}} Z_{b}^{\circ},
                \end{tikzcd}
                $$
                where the right vertical map is induced by the morphisms $\phi_{j} \colon \widetilde{Y}_{j} \to Y_{j}$ and $\varphi_{k} \colon \widetilde{Z}_{k} \to Z_{k}$ over $Y$ such that for each point $y \in Y^{\circ}$, the morphisms $\phi_{j}$ and $\varphi_{k}$, restricted to the fibers over $y$, give resolutions of singularities.
                \item $\mu_{1}$ is a birational morphism, $h$ is a generically finite morphism, and $\widetilde{f} \circ h$ satisfies the standard normal crossing assumptions.
                \item The variations of Hodge structures $R^{\dim p}p_{\ast} \CC_{A^{\circ}}$, $R^{\dim \widetilde{q}_{j}}\widetilde{q}_{j\ast} \CC_{\widetilde{Y}_{j}^{\circ}}$, $R^{\dim \widetilde{r}_{k}}\widetilde{r}_{k\ast} \CC_{\widetilde{Z}_{k}^{\circ}}$, $R^{\dim \widetilde{f}} \widetilde{f}\lsta \CC_{\widetilde{X}^{\circ}}$, and $R^{\dim \widetilde{f}}(\widetilde{f} \circ h)_{\ast} \CC_{{X'''}^{\circ}}$ have unipotent local monodromies along $\Sigma$.
            \end{enumerate}
            Note that $M_{X/ Y} = M_{\widetilde{X}/ Y} = M_{X'''/Y} = M_{X''/Y}$ by Proposition \ref{prop:Fujino-cover-formula} and Remark \ref{rema:moduli-upto-res-of-X}. Note that
            $$M_{X''/Y} = M_{A/ Y} + \sum_{j=1}^{a} M_{\widetilde{Y}_{j}/ Y} + \sum_{k=1}^{b} M_{\widetilde{Z}_{k}/ Y} = M_{A/ Y} + \sum_{j=1}^{a} M_{Y_{j}/ Y} + \sum_{k=1}^{b} M_{Z_{k}/ Y}. $$
            We already know that $M_{A/ Y}$ and $M_{Y_{j}/ Y}$ are semi-ample. Therefore, Conjecture \ref{conj:semi-ample-CY} (2) in relative dimension $\leq n$ implies Conjecture \ref{conj:semi-ample-allK-triv} (2) in relative dimension $n$.
        \end{proof}

        \begin{rema}
            Note that we used wider classes of varieties when we passed from the properties of geometric generic fiber to those of the general fiber. This is because the definitions of IHS and CY in \cite{Greb-Guenancia-Kebekus:kltBeauville-Bog}*{Definition 1.3} have a condition on the reflexive Hodge numbers for `every' étale in codimension 1 cover, which involves a priori infinite data. It is also not clear to use the description of IHS and CY using holonomy groups in the sense of \cite{Greb-Guenancia-Kebekus:kltBeauville-Bog}*{Proposition 12.10} since one should choose an abstract field isomorphism $\CC \simeq \overline{k(\eta)}$ in order to perform such a decomposition in the proof of Theorem \ref{theo:reduction-to-CY}. It would be interesting if one could overcome this technical point.
        \end{rema}

        {\bf Acknowledgements.} The author would like to thank his advisor Mircea Musta\c{t}\u{a} for support and helpful discussions. The author would also like to thank Osamu Fujino and Stefano Filipazzi for useful comments and answering his questions.
	    
        \bibliographystyle{alpha}
        \bibliography{Reference}

\end{document}